\documentclass{amsart}

\usepackage{amsmath,amssymb}
\usepackage{amsthm}
\usepackage[alphabetic]{amsrefs}
\usepackage{indentfirst}
\usepackage{mathrsfs}
\usepackage{graphicx}
\usepackage{hyperref}
\usepackage{xcolor}
\usepackage[margin=1.5in]{geometry}
\usepackage{graphicx,tikz}

\newtheorem*{thm}{Theorem}

\newtheorem*{lemma}{Lemma}

\theoremstyle{definition}

\theoremstyle{remark}

\DeclareMathOperator{\diam}{diam}

\begin{document}

\title[]{ Leaky roots and stable Gauss-Lucas Theorems}
\keywords{Gauss-Lucas theorem, roots, zeroes, derivatives, critical points.}
\subjclass[2010]{26C10, 30C15.}

\author[]{Trevor J. Richards}
\address{Department of Mathematics and Computer Science, Monmouth College, Monmouth,
IL 61462, USA.}
\email{trichards@monmouthcollege.edu}

\author[]{Stefan Steinerberger}
\address{Department of Mathematics, Yale University, New Haven, CT 06511, USA}
\email{stefan.steinerberger@yale.edu}
\thanks{S.S. is supported by the NSF (DMS-1763179) and the Alfred P. Sloan Foundation.}

\begin{abstract}
Let $p:\mathbb{C} \rightarrow \mathbb{C}$ be a polynomial. The Gauss-Lucas theorem states that its critical points, $p'(z) = 0$, are contained in the convex hull of its roots. A recent quantitative
version Totik shows that if almost all roots are contained in a bounded convex domain $K \subset \mathbb{C}$, then almost all roots of the derivative $p'$ are 
in a $\varepsilon-$neighborhood $K_{\varepsilon}$ (in a precise sense). We prove another quantitative version: if a polynomial $p$ has $n$ roots in $K$ and $\lesssim c_{K, \varepsilon} (n/\log{n})$ roots outside of $K$,
then $p'$ has at least $n-1$ roots in $K_{\varepsilon}$. This establishes, up to a logarithm, a conjecture of the first author: we also discuss an open problem whose solution would imply the full conjecture.
\end{abstract}

\maketitle

\section{Introduction and result}
\subsection{Introduction.} The Gauss-Lucas Theorem, stated by Gauss \cite{gauss} in 1836 and proved by Lucas \cite{lucas} in 1879, says that if $p_n:\mathbb{C} \rightarrow \mathbb{C}$ is a polynomial of degree $n$, then the $n-1$ zeroes of $p_n'$ lie inside the convex hull of the $n$ zeros of $p_n$. This has been of continued interest \cite{aziz, bray, bruj, bruj2, branko, dimitrov, joyal, kalman, malamud, marden, pawlowski, pereira, rav, trevor, siebeck, schm, specht, stein, walsh}.  Totik \cite{totik} recently showed that, for sequences of polynomials $p_n$ with $\deg(p_n) \rightarrow \infty$, that if $(1-o(1))\deg{p_n}$ roots of $p_n$ lie inside a convex domain $K$, then any fixed $\varepsilon-$neighborhood of $K$ contains $(1-o(1))\deg{p_n}$ roots of $p'$. 

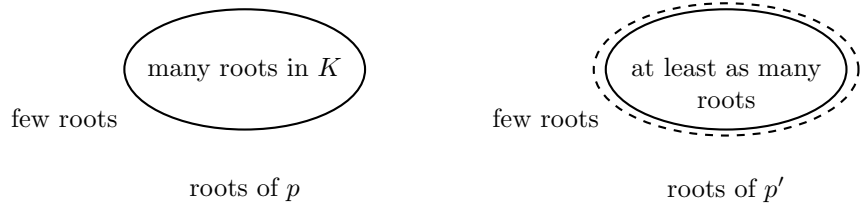
\begin{figure}[h!]
\begin{center}
\begin{tikzpicture}[scale=0.8]
\draw [ thick] (0,0) ellipse (2cm and 1cm);
\node at (0,0) {many roots in $K$};
\node at (-3,-0.8) {few roots};
\node at (0,-2) {roots of $p$};
\draw [ thick] (8,0) ellipse (2cm and 1cm);
\draw [ thick, dashed] (8,0) ellipse (2.2cm and 1.1cm);
\node at (8,0) {at least as many};
\node at (8,-0.5) {roots};
\node at (8-3,-0.8) {few roots};
\node at (8,-2) {roots of $p'$};
\end{tikzpicture}
\end{center}
\vspace{-10pt}
\caption{The structure of asymptotic Gauss-Lucas theorems.}
\end{figure}
These results rely on potential theory and are asymptotic in nature. It would be
of interest to have stable non-asymptotic statements showing that the roots of $p_n'$ cannot move very much outside of $K$, i.e. statements along the lines of
$$ \# \left\{z \in K_{\varepsilon}: p_n'(z) = 0 \right\} \geq \# \left\{z \in K_{}: p_n(z) = 0 \right\} -1,$$
where $K_{\varepsilon}$ is the $\varepsilon-$neighborhood of $K$. A schematic picture of such results is given in Fig. 1.
The first author has proposed a nonasymptotic formulation in \cite{trevor}.
\begin{quote}
\textbf{Conjecture} (Richards \cite{trevor})\textbf{.} For every $\varepsilon > 0$ there exists a constant $c_{K, \varepsilon} > 0$ depending only on the bounded convex domain $K \subset \mathbb{C}$ and $\varepsilon>0$ such that for all polynomials $p_n:\mathbb{C} \rightarrow \mathbb{C}$ of degree $n$ the following holds: if $p_n$ has at most $c_{K, \varepsilon}n$ roots outside of $K$, then 
$$ \# \left\{z \in K_{\varepsilon}: p_n'(z) = 0 \right\} \geq \# \left\{z \in K_{}: p_n(z) = 0 \right\} - 1.$$
\end{quote}

 The first author \cite{trevor} proved the conjecture if the number of roots outside of $K$ is $\lesssim_{K, \varepsilon} \sqrt{n}$.
The purpose of our paper is to prove the conjecture up to a logarithm and to discuss a new conjecture about the geometry of Coulomb potentials in $\mathbb{C}$ that seems of intrinsic interest and would imply
the full conjecture.

\begin{thm} Let $K \subset \mathbb{C}$ be a bounded convex domain. For every $\varepsilon > 0$ there exists a constant $c_{K, \varepsilon} > 0$ such that for all polynomials $p_n: \mathbb{C} \rightarrow \mathbb{C}$
of degree $n$ the following holds: if 
$$ \# \left\{z \in \mathbb{C} \setminus K: p_n(z) = 0 \right\} \leq c_{K, \varepsilon} \frac{n}{\log{n}},$$
then 
$$ \# \left\{z \in K_{\varepsilon}: p_n'(z) = 0 \right\} \geq \# \left\{z \in K_{}: p_n(z) = 0 \right\} -1.$$
\end{thm}
The full conjecture would follow from the truth of the following curious conjecture that we believe to be of interest in itself: let $\gamma:[0,1] \rightarrow \mathbb{C}$ be a curve in
the complex plane normalized to $\gamma(0) = 0$ and $\gamma(1) = 1$, let $n \in \mathbb{N}$ be a parameter and let $\left\{z_1, \dots, z_m\right\} \subset \mathbb{C} \setminus \gamma([0,1])$ be a set
of $m$ charges in the complex plane. We conjecture~\footnote{The authors have privately received a proof of a generalization of this conjecture from Vilmos Totik, who showed that the same conclusion holds for general measures (not just sums of point masses).  As noted above, this in turn implies the conjectured approximate Gauss--Lucas theorem which originated in~\cite{trevor}.  We look forward to seeing Totik's general result in print in the future.} that $m$ charges cannot `supercharge' a curve. More precisely, we conjecture that if
$$ \min_{0 \leq t \leq 1}{ \left| \sum_{\ell=1}^{m}{\frac{1}{\gamma(t)-z_{\ell}}} \right| } \geq n \qquad \mbox{then} \qquad m \gtrsim n.$$
This, if true, would then imply the conjecture via our approach outlined below.
We prove $m \gtrsim n/\log{n}$ (under slightly weaker assumptions for which this scaling is sharp).
\begin{figure}[h!]
\begin{center}
\begin{tikzpicture}[scale = 2]
\filldraw (0,0) circle (0.04cm);
\node at (0,-0.3) {0};
\draw [thick] (0,0) to[out=30, in=120] (2,0);
\filldraw (2,0) circle (0.04cm);
\node at (2,-0.3) {1};
\filldraw (0.2,0) circle (0.03cm);
\filldraw (1,0.5) circle (0.03cm);
\filldraw (1.7,0.38) circle (0.03cm);
\end{tikzpicture}
\end{center}
\caption{Strategically placed point charges charging a curve.}
\end{figure}
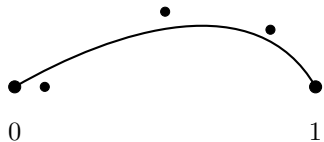

One could also ask the same question with the $m$ point charges replaced by a measure $\mu$. 

\section{A discrete Lemma}
Let $\gamma:[0,1] \rightarrow \mathbb{C}$ be a curve whose endpoints satisfy
$|\gamma(0) - \gamma(1)| = 1$. Let $\left\{z_1, \dots, z_m \right\} \subset \mathbb{C} \setminus \gamma([0,1])$ be arbitrary. We prove that a finite number of points cannot
strongly charge the curve via convolution with a $|x|^{-1}$ potential unless there is a sufficient number of them.
\begin{lemma} We have
$$ \min_{0 \leq t \leq 1}{ \sum_{\ell=1}^{m}{\frac{1}{|\gamma(t)-z_{\ell}|}  }} \leq 60 m \log{m}$$
\end{lemma}
It is easy to see that the Lemma is sharp up to the constant: let us pick $\gamma(t) = t + 0i$ for $0 \leq t \leq 1$ and $z_j = j/m + i/m$. Then
$$ \min_{0 \leq t \leq 1}{ \sum_{\ell=1}^{m}{\frac{1}{|\gamma(t)-z_{\ell}|}  }} \gtrsim \sum_{j=1}^{m}{ \frac{m}{j}} \sim m \log{m}.$$

\begin{proof}[Proof of Lemma 1] We will bound the quantity from above by restricting the curve
$\gamma$ to its $x-$coordinate and the complex number to its real part to conclude that
$$  \min_{0 \leq t \leq 1}{ \sum_{\ell=1}^{m}{\frac{1}{|\gamma(t)-z_{\ell}|}  }} \leq  \min_{0 \leq t \leq 1}{ \sum_{\ell=1}^{m}{\frac{1}{|\Re \gamma(t)- \Re z_{\ell}|}  }} \leq   \min_{0 \leq t \leq 1}{ \sum_{\ell=1}^{m}{\frac{1}{|t - \Re z_{\ell}|}  }} .$$
By `tying the ends of the interval together', we can replace that problem by yet another problem: given $\left\{x_1, \dots, x_m\right\} \subset \mathbb{T}$ (where $\mathbb{T}$ is the one-dimensional torus with length 1), our desired result would follow from knowing that 
$$\min_{x \in \mathbb{T}}{ \sum_{\ell=1}^{m}{\frac{1}{|x - x_{\ell}|}  }} \leq 60 m \log{m},$$
where $\left| \cdot \right|$ denotes the toroidal distance $|x| = \min\left\{x, 1-x\right\}$. We re-interpret the problem yet again: by writing it as a convolution and using $\delta_x$ to denote the Dirac measure in a point $x$, our desired result would follow from showing that
$$ \min_{y \in \mathbb{T}} \left( \frac{1}{|x|} *\left( \sum_{\ell = 1}^{m}{\delta_{x_{\ell}}}\right) \right)(y) \leq 60 m \log{m}.$$
We now replace the singular kernel by a slightly less singular function 
$$ f_m(x) = \begin{cases} |x|^{-1} \qquad &\mbox{if}~x \geq 1/(20m) \\ 20m \qquad &\mbox{otherwise.} \end{cases}$$
We observe that this function is integrable and
\begin{align*}
 \int_{\mathbb{T}} f_m * \left( \sum_{\ell = 1}^{m}{\delta_{x_{\ell}}} \right)  dx &= m \int_{\mathbb{T}}{f_m(x) dx} \\
&\leq m + 2 m  \int_{(20m)^{-1}}^{1/2}{\frac{dx}{x}} \leq  2 m  \int_{(20m)^{-1}}^{1}{\frac{dx}{x}} = 2 m \log{(20m)}.
\end{align*}

Markov's inequality now implies that the set where the convolution is large is small
$$ \left| \left\{ x \in \mathbb{T}: \sum_{\ell=1}^{m}{f_m(x-x_{\ell})} \geq 20 m \log{(20m)} \right\} \right| \leq \frac{1}{10}.$$
However, it is also easy to see that the set of points really close to one of the $m$ points cannot be too large
$$ \left| \left\{x \in \mathbb{T}: \min_{1 \leq \ell \leq m}{|x-x_{\ell}|} \leq \frac{1}{10m} \right\} \right| \leq \frac{1}{5}$$
and therefore there exists a point $y \in \mathbb{T}$ that is at distance at least $(10m)^{-1}$ from $\left\{x_1, \dots, x_m\right\}$ and for which
$$\sum_{\ell=1}^{m}{f_m(y-x_{\ell})} \leq 20 m \log{(20m)} \leq 60m\log{m}.$$
However, for points at that distance, $f_m$ and $|x|^{-1}$ coincide and this proves the desired claim.
\end{proof}

\section{Proof of the Theorem}

\begin{proof} We will assume that polynomial $p$ has $n+m$ roots where the $n$ roots in $K$ are given by $a_1, \dots, a_n \in K$ and that the roots $a_{n+1}, \dots, a_{n+m}$ are contained outside of $K$.
The roots of a polynomial depend smoothly on the coefficients and, conversely, the coefficients of a polynomial depend smoothly on the roots. This allows us to assume for simplicity that all roots are distinct and
that $p'$ and $p$ have no common roots.
We note that roots of the derivative have to satisfy
$$ \frac{p'(z)}{p(z)} = \sum_{k=1}^{n+m}{\frac{1}{z-a_k}} = 0.$$
We write the polynomial as $p(z) = q(z) r(z)$ where
$$ q(z) = \prod_{k=1}^{n}{(z-a_k)} \qquad \mbox{and} \qquad r(z) = \prod_{k=n+1}^{n+m}{(z-a_k)}.$$
Then $p'(z) = q'(z)r(z) + q(z) r'(z)$. We will now introduce the set, for any $\delta >0$,
$$ A_{\delta} = \left\{ z \in \mathbb{C}: \left| \frac{q'}{q} \right| \leq \left| \frac{r'}{r} \right|  + \frac{\delta}{|r|} \right\}.$$
We observe that every root of $p'$ is contained in $A_{\delta}$ for all $\delta > 0$: this is because if $z$ is a root of $p'$,
then $q'(z)/q(z) = - r'(z)/r(z)$. Simple far-field asymptotics also show that $A_{\delta}$ is bounded (for $\delta$ sufficiently small if $m=1$ and unconditionally if $m \geq 2$ as long as $m < n$).
If $A$ is a connected component of $A_{\delta}$, then 
$$ |q'(z) r(z) | =  |q(z) r'(z)|  + \delta |q(z)| \qquad \mbox{and thus} \qquad  |q'(z) r(z) |  > |q(z) r'(z)| \quad \mbox{on}~\partial A.$$
Rouch\'{e}'s theorem implies that the number of roots of $p'(z) = q'(z)r(z) + r(z) q'(z)$ in $A$ is the same as the number of roots of $q'(z) r(z)$ in $A$. 
This has an important implication: it identifies regions in the complex plane, the connected components of $A_{\delta}$ (for sufficiently small $\delta$), where the number of roots of $p'$ and the number of the roots of $q'(z)r(z)$ exactly coincide. 
We are interested in the number of roots of $p'$ in $\mathbb{C} \setminus K_{\varepsilon}$ and want to show that this number is at most $m$. Let now $z_0 \in \mathbb{C} \setminus K_{\varepsilon}$, let $p'(z_0) = 0$ and let $A$ be the connected component of $A_{\delta}$ containing $z_0$. If $A$ does not intersect $K$, then $q'$ never vanishes in $A$ (this follows from the Gauss-Lucas theorem) and every root of
$p'$ corresponds to a root of $r$ (of which there are $m$ in total). It thus suffices to establish a condition under which no connected component of $A_{\delta}$ intersects both $\mathbb{C} \setminus K_{\varepsilon}$ and $K$.
We will now establish such a condition.

\begin{figure}[h!]
\begin{center}
\begin{tikzpicture}[scale=1]
\draw [thick] (0,0) ellipse (2cm and 1cm);
\node at (0,0) {$K$};
\draw [thick] (0,0) ellipse (2.8cm and 1.8cm);
\node at (2,-0.8) {$K_{\varepsilon}$};
\draw [thick, dashed] (-2,-2) to[out=300, in = 270] (-1, -2) to[out=90, in = 0] ( -1.5, -1) to[out=180, in=120] (-2,-2);
\draw [thick, dashed] (-2-2,-2+2) to[out=300, in = 270] (-1-2, -2+2) to[out=90, in = 0] ( -1.5-2, -1+2) to[out=180, in=120] (-2-2,-2+2);
\draw [thick, dashed] (3-2,-2+2) to[out=300, in = 270] (3-2, -2+2) to[out=90, in = 0] (3-2, -1+2) to[out=180, in=120] (3-2,-2+2);
\draw [thick, dashed] (4, 1) ellipse (0.5cm and 1cm);
\end{tikzpicture}
\end{center}
\caption{$K$, the $\varepsilon-$neighborhood $K_{\varepsilon}$ and the set $A_{\delta}$ (dashed). Every connected component of $A_{\delta}$ contains exactly as many roots of $p'$ as it has roots of $r$.}
\end{figure}
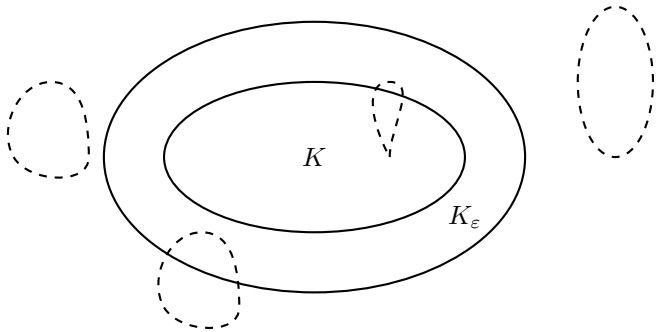

Suppose the statement fails and the connected component $A$ of $A_{\delta}$ intersects both $\mathbb{C} \setminus K_{\varepsilon}$ and $K$. It then contains a curve connecting $\mathbb{C} \setminus K_{\varepsilon}$ and $K$ (a curve that is at least of length $\varepsilon$) for all $\delta >0$. We will now construct a lower bound on
$$ \left| \frac{q'(z)}{q(z)} \right| = \left|  \sum_{k=1}^{n}{ \frac{1}{ z - a_k}} \right| \qquad \mbox{for}~z \in \mathbb{C} \setminus K.$$
By rotational and translational variance, we can again assume that $z \in \mathbb{R}$, that $(0,0) \in K$ is the closest point to $z$ in $K$, and 
$$z> \sup_{y \in K}{\mbox{Re}~y} = 0.$$

We now estimate $|q'/q|$ at $z$ by $$\left|\sum_{k=1}^n\dfrac{1}{z-a_k}\right|\geq\displaystyle\Re\left(\sum_{k=1}^n\dfrac{1}{z-a_k}\right)=\sum_{k=1}^n\dfrac{\Re(z-\overline{a_k})}{|z-a_k|^2}.$$  It now follows from the normalization described above that $$\left|\sum_{k=1}^n\dfrac{1}{z-a_k}\right|\geq\sum_{k=1}^n\dfrac{d(z,K)}{(d(z,K)+\operatorname{diam}(K))^2}=\dfrac{n\cdot d(z,K)}{(d(z,K)+\operatorname{diam}(K))^2}.$$

The assumption of $A$ intersecting both $\mathbb{C}\setminus K_\epsilon$ and $K$ implies that $A$ contains a path $\gamma$ whose endpoints are distance $\geq\varepsilon/2$ from each other and which lies a distance at least $\varepsilon/2$ from $K$.  The lower bound on $|q'(z)/q(z)|$ obtained above, along with the definition of the set $A$ thus gives us that for $z$ in $\gamma$, $$\left|\dfrac{p'(z)}{p(z)}\right|=\left|\displaystyle\sum_{k=n+1}^{n+m}\dfrac{1}{z-a_k}\right|\geq\frac{1}{2}\dfrac{n\varepsilon}{(\varepsilon+\operatorname{diam}(K))^2},$$ and thus for $\varepsilon<\operatorname{diam}(K)$ and some constant $c_K>0$, $$\left|\displaystyle\sum_{k=n+1}^{n+m}\dfrac{1}{z-a_k}\right|>c_Kn\varepsilon\ \ \ \text{ for $z\in\gamma$}.$$

Rescaling by a factor of $\varepsilon^{-1}$ and using the Lemma implies that
 $$c_Kn\varepsilon\leq\dfrac{60m\log m}{\varepsilon} \qquad \text{forcing} \qquad m\gtrsim_{K,\varepsilon}\dfrac{n}{\log n}.$$

\end{proof}

\end{document}